\documentclass[10pt]{article}
\usepackage{latexsym,amsmath,amsfonts,graphics}
\newcommand{\ignore}[1]{}

\begin{document}
\begin{center}
{\LARGE\textbf{Bounds of distance Estrada index of graphs}}\\
\bigskip
\bigskip
Yilun Shang\\
Department of Mathematics\\
Tongji University, Shanghai 200092, China\\
e-mail: \texttt{shyl@tongji.edu.cn}
\end{center}

\smallskip
\begin{abstract}

Let $\lambda_1,\lambda_2,\cdots,\lambda_n$ be the eigenvalues of the
distance matrix of a connected graph $G$. The distance Estrada index
of $G$ is defined as $DEE(G)=\sum_{i=1}^ne^{\lambda_i}$. In this
note, we present new lower and upper bounds for $DEE(G)$. In
addition, a Nordhaus-Gaddum type inequality for $DEE(G)$ is given.

\bigskip



\end{abstract}

\bigskip
\normalsize

\section{Introduction}

Let $G=(V(G),E(G))$ be a simple connected graph with vertex set
$V(G)=\{v_1,v_2,\cdots,v_n\}$ and edge set $E(G)$. The adjacency
matrix of $G$ is $A(G)=(a_{ij})\in\mathbb{R}^{n\times n}$ (or $A$
for short), where $a_{ij}=1$ if two vertices $v_i$ and $v_j$ are
adjacent in $G$ and $a_{ij}=0$ otherwise. The eigenvalues of $A$ are
real, and can be ordered as
$\lambda_1(A)\ge\lambda_2(A)\ge\cdots\ge\lambda_n(A)$. The distance
matrix of $G$ is a symmetric matrix
$D(G)=(d_{ij})\in\mathbb{R}^{n\times n}$ (or $D$ for short) in which
$d_{ij}$ denotes the length of shortest path between two vertices
$v_i$ and $v_j$. Since $D$ is real symmetric, its eigenvalues
(called distance eigenvalues or distance spectra) can also be
arranged in non-increasing order as
$\lambda_1(D)\ge\lambda_2(D)\ge\cdots\ge\lambda_n(D)$. We refer the
reader to \cite{1} for a comprehensive survey on distance
eigenvalues.

The Estrada index of graph $G$, put forward by Estrada \cite{2}, is
defined as
$$
EE(G)=\sum_{i=1}^ne^{\lambda_i(A)}.
$$
This graph-spectrum-based invariant has found a number of
applications in chemistry, physics, and complex networks. For
example, it is used as a measure for the degree of folding of long
chain polymeric molecules \cite{3}. It also characterizes the
centrality \cite{4} as well as robustness \cite{5} of complex
networks. For various mathematical properties of the Estrada index,
see e.g. \cite{8,6,7,a}.

Quite recently, in full analogy with the Estrada index, the distance
Estrada index of a connected graph $G$ was introduced in \cite{9} as
\begin{equation}
DEE(G)=\sum_{i=1}^ne^{\lambda_i(D)}.\label{1}
\end{equation}
Let $K_n$ be the complete graph on $n$ vertices. Some upper and
lower bounds for $DEE(G)$ were established as follows.

\smallskip
\noindent\textbf{Theorem 1.} \cite{9}\itshape \quad Let $G$ be a
connected graph with $n$ vertices and $m$ edges, and $\rho=\rho(G)$
the diameter of $G$. Then
\begin{equation}
\sqrt{n^2+4m}\le DEE(G)\le n-1+e^{\rho\sqrt{n(n-1)}}.\label{2}
\end{equation}
Equality holds on both sides of (\ref{2}) if and only if $G=K_1$.
\normalfont
\smallskip

\smallskip
\noindent\textbf{Theorem 2.} \cite{10}\itshape \quad Let $G$ be a
connected graph with $n\ge2$ vertices and $m$ edges. Then
\begin{equation}
DEE(G)\ge e^{2(n-1)-\frac{2m}{n}}+e^{-2(n-1)+\frac{2m}{n}}+n-2
\label{3}
\end{equation}
with equality if and only if $G=K_2$. \normalfont
\smallskip

Some results relating $DEE(G)$ to the Winer index and graph energy
can be found in \cite{10,9}. Moreover, the distance Estrada index
for strongly quotient graphs and Erd\H{o}s-R\'enyi random graphs
were discussed in \cite{11} and \cite{12}, respectively. In this
note, we establish some new bounds for $DEE(G)$ involving diameter,
maximum degree, and second maximum degree. Our bounds improve some
results in \cite{10,9}. Furthermore, a Nordhaus-Gaddum type
inequality for $DEE(G)$ is presented.

\section{Some lemmas}

We list some useful lemmas for distance spectra in this section.

\smallskip
\noindent\textbf{Lemma 1.} \cite{15}\itshape \quad Let $G$ be a
connected graph with $n$ vertices and $m$ edges. Then
$\sum_{i=1}^n\lambda_i(D)=0$ and
$\sum_{i=1}^n\lambda_i^2(D)=2\sum_{i<j}d_{ij}^2$. \normalfont
\smallskip

The next result relies on a special relation between the adjacency
and distance matrices of graphs having diameter 2. It has been
applied in deriving $DEE(G)$ in dense random graph settings
\cite{12}.

\smallskip
\noindent\textbf{Lemma 2.} \cite{13,14}\itshape \quad Let $G$ be an
$r$-regular graph on $n$ vertices with diameter at most 2 and
adjacency eigenvalues $\lambda_1(A)=r$,
$\lambda_2(A),\lambda_3(A),\cdots,\lambda_n(A)$. Then the distance
eigenvalues of $G$ are $2n-2-r$, $-2-\lambda_2(A)$,
$-2-\lambda_3(A)$, $\cdots, -2-\lambda_n(A)$. \normalfont
\smallskip

\smallskip
\noindent\textbf{Lemma 3.} \cite{16}\itshape \quad Let $G$ be a
connected graph on $n$ vertices with maximum degree $\Delta_1$ and
second maximum degree $\Delta_2$. Then
$$
\lambda_1(D)\ge\sqrt{(2n-2-\Delta_1)(2n-2-\Delta_2)}
$$
with equality if and only if $G$ is a regular graph with diameter at
most 2. \normalfont
\smallskip

More details on extremal values for the largest distance eigenvalue
of a graph can be found in e.g. \cite{18,19,16}. Let
$K_{n_1,n_2,\cdots,n_s}$ denote the complete $s$-partite graph. The
following result characterizes the least distance eigenvalue.

\smallskip
\noindent\textbf{Lemma 4.} \cite{17}\itshape \quad Let $G$ be a
graph on $n$ vertices. For $n\ge2$, $\lambda_n(D)=-1$ if and only if
$G=K_n$. For $n\ge3$, $\lambda_n(D)=-2$ if and only if
$G=K_{n_1,n_2\cdots,n_s}$ for some $s\in[2,n-1]$ with
$\sum_{i=1}^sn_i=n$. Moreover, for $n\ge3$, if $G\not=K_n$ and
$G\not=K_{n_1,n_2\cdots,n_s}$ for some $s\in[2,n-1]$ with
$\sum_{i=1}^sn_i=n$, then $\lambda_n(D)<-2.383$. \normalfont
\smallskip

\section{Bounds for distance Estrada index}

\smallskip
\noindent\textbf{Theorem 3.} \itshape \quad Let $G$ be a connected
graph on $n\ge2$ vertices with maximum degree $\Delta_1$ and second
maximum degree $\Delta_2$. Then
\begin{equation}
DEE(G)\ge
e^{\sqrt{(2n-2-\Delta_1)(2n-2-\Delta_2)}}+(n-1)e^{-\sqrt{\big(2-\frac{\Delta_1}{n-1}\big)\big(2-\frac{\Delta_2}{n-1}\big)}}\label{4}
\end{equation}
with equality if and only if $G=K_n$. \normalfont
\smallskip

\noindent\textbf{Proof.} Using Lemma 1 and the arithmetic-geometric
inequality, we obtain
\begin{eqnarray*}
DEE(G)&=&e^{\lambda_1(D)}+\sum_{i=2}^ne^{\lambda_i(D)}\\
&\ge&e^{\lambda_1(D)}+(n-1)\Big(e^{\sum_{i=2}^n\lambda_i(D)}\Big)^{\frac{1}{n-1}}\\
&=&e^{\lambda_1(D)}+(n-1)e^{-\frac{\lambda_1(D)}{n-1}},
\end{eqnarray*}
with equality if and only if
$\lambda_2(D)=\lambda_3(D)=\cdots=\lambda_n(D)$.

For $x\ge0$, define $f(x)=e^x+(n-1)e^{-\frac{x}{n-1}}$. It is easy
to see that $f'(x)=e^x-e^{-\frac{x}{n-1}}\ge0$ for any $x\ge0$. By
Lemma 3, we have
$\lambda_1(D)\ge\sqrt{(2n-2-\Delta_1)(2n-2-\Delta_2)}\ge0$.
Therefore,
$$
DEE(G)\ge f(\lambda_1(D))\ge
f\big(\sqrt{(2n-2-\Delta_1)(2n-2-\Delta_2)}\big),
$$
and (\ref{4}) follows.

To see the sharpness of (\ref{4}), note that the complete graph
$K_n$ has distance spectrum $\{n-1,-1,-1,\cdots,-1\}$. Hence, if
$G=K_n$, we obtain $DEE(G)=e^{n-1}+(n-1)e^{-1}$, and the equality
holds in (\ref{4}).

Conversely, if the equality holds in (\ref{4}) then $G$ must be a
regular graph, say, $r$-regular, with diameter at most 2 by Lemma 3.
Employing Lemma 2 and the fact that
$\lambda_2(D)=\lambda_3(D)=\cdots=\lambda_n(D)$, we conclude that
$\lambda_1(A)=r$, $\lambda_2(A)=\lambda_3(A)=\cdots=\lambda_n(A)$.
Since the trace of $A$ is zero, we know that
$\lambda_1(A)>\lambda_2(A)$. It is well known that a connected graph
with two distinct adjacency eigenvalues must be complete, which
completes the proof. $\Box$

It is evident that our bound is better than the lower bound in
(\ref{2}). Notice that our bound is incomparable to the bound in
(\ref{3}). Nevertheless, there are more graphs which attain our
bound.

In 1956, Nordhaus and Gaddum \cite{20} presented lower and upper
bounds on the sum of the chromatic number of a graph and its
complement, in terms of the order of the graph. Here, we give a
Nordhaus-Gaddum type result for the distance Estrada index.

\smallskip
\noindent\textbf{Theorem 4.} \itshape \quad Let $G$ be a connected
graph on $n\ge2$ vertices with a connected complement
$\overline{G}$. Then
\begin{equation}
DEE(G)+DEE(\overline{G})>2e^{\frac{3(n-1)}{2}}+2e^{-\frac{3(n-1)}{2}}+2n-4.\label{5}
\end{equation}
\normalfont
\smallskip

\noindent\textbf{Proof.} Define a function
$g(n,m)=e^{2(n-1)-\frac{2m}{n}}+e^{-2(n-1)+\frac{2m}{n}}+n$.

Suppose that $|E(G)|=m$. Hence, $|E(\overline{G})|={n\choose2}-m$.
Since $\overline{G}$ is connected, it follows from Theorem 2 that
\begin{eqnarray*}
DEE(G)+DEE(\overline{G})&>&g(n,m)+g\left(n,{n\choose2}-m\right)-4\\
&=&e^{2(n-1)-\frac{2m}{n}}+e^{-2(n-1)+\frac{2m}{n}}+e^{2(n-1)-\frac{n(n-1)-2m}{n}}\\
&&+e^{-2(n-1)+\frac{n(n-1)-2m}{n}}+2n-4.
\end{eqnarray*}

We have $e^{x_1}+e^{x_2}\ge2e^{\frac{x_1+x_2}{2}}$ since $e^x$ is
convex. Therefore, we obtain
$$
DEE(G)+DEE(\overline{G})>2e^{\frac{3(n-1)}{2}}+2e^{-\frac{3(n-1)}{2}}+2n-4
$$
as desired. $\Box$

\smallskip
\noindent\textbf{Theorem 5.} \itshape \quad Let $G$ be a connected
graph with $n\ge2$ vertices, and $\rho=\rho(G)$ the diameter of $G$.
Then
\begin{equation}
DEE(G)<n-1+e^{\sqrt{n(n-1)\rho^2-1}}.\label{6}
\end{equation}
\normalfont
\smallskip

\noindent\textbf{Proof.} Denote by $n_+$ the number of positive
eigenvalues of $D$. We obtain
\begin{eqnarray}
DEE(G)=\sum_{i=1}^ne^{\lambda_i(D)}&<&n-n_++\sum_{i=1}^{n_+}e^{\lambda_i(D)}\nonumber\\
&=&n-n_++\sum_{i=1}^{n_+}\sum_{k\ge0}\frac{\lambda_i^k(D)}{k!}\nonumber\\
&=&n+\sum_{k\ge1}\frac{1}{k!}\sum_{i=1}^{n_+}\lambda_i^k(D)\nonumber\\
&=&n+\sum_{k\ge1}\frac{1}{k!}\sum_{i=1}^{n_+}(\lambda_i^2(D))^{\frac{k}{2}}\nonumber\\
&\le&n+\sum_{k\ge1}\frac{1}{k!}\bigg(\sum_{i=1}^{n_+}\lambda_i^2(D)\bigg)^{\frac{k}{2}},\label{7}
\end{eqnarray}
where the first inequality holds due to the strict monotonicity of
$e^x$.

In view of Lemma 1 and (\ref{7}), we deduce
$$
DEE(G)<n+\sum_{k\ge1}\frac{1}{k!}\bigg(2\sum_{i<j}d_{ij}^2-\sum_{i=n_++1}^n\lambda_i^2(D)\bigg)^{\frac{k}{2}}.
$$
Since $G$ is a connected graph on $n\ge2$ vertices, Lemma 4 implies
that $\lambda_n(D)\le-1$. Accordingly,
\begin{eqnarray*}
DEE(G)&<&n+\sum_{k\ge1}\frac{1}{k!}\bigg(2\sum_{i<j}d_{ij}^2-1\bigg)^{\frac{k}{2}}\\
&\le&n+\sum_{k\ge1}\frac{1}{k!}\big(n(n-1)\rho^2-1\big)^{\frac{k}{2}}\\
&=&n-1+e^{\sqrt{n(n-1)\rho^2-1}},
\end{eqnarray*}
where the second inequality holds since $d_{ij}\le\rho$, and the
last equality is because of the power-series expansion of $e^x$. The
proof is complete. $\Box$

Obviously, our upper bound is better than that in (\ref{2}).

To conclude this note, we mention a result relating $DEE(G)$ to
$EE(\overline{G})$ for regular graphs, which is a direct corollary
of Lemma 2.

\smallskip
\noindent\textbf{Theorem 6.} \itshape \quad Let $G$ be an
$r$-regular graph on $n$ vertices with diameter at most 2. Then
\begin{equation}
DEE(G)=e^{2n-r-2}-e^{n-r-2}+e^{-1}EE(\overline{G}).\label{8}
\end{equation}
\normalfont
\smallskip

\noindent\textbf{Proof.} It is well known that (see e.g.
\cite[p.172]{21}) the adjacency eigenvalues of $\overline{G}$ are
$\{n-r-1,-1-\lambda_2(A(G)),-1-\lambda_3(A(G)),\cdots,-1-\lambda_n(A(G))\}$.

Hence,
$EE(\overline{G})=\sum_{i=1}^ne^{\lambda_i(A(\overline{G}))}=e^{n-r-1}+e(e^{-2-\lambda_2(A(G))}+\cdots+e^{-2-\lambda_n(A(G))})$.
Thanks to Lemma 2, we have
\begin{eqnarray*}
DEE(G)&=&\sum_{i=1}^ne^{\lambda_i(D)}=e^{2n-r-2}+e^{-2-\lambda_2(A(G))}+\cdots+e^{-2-\lambda_n(A(G))}\\
&=&e^{2n-r-2}+e^{-1}\big(EE(\overline{G})-e^{n-r-1}\big),
\end{eqnarray*}
which readily implies the desired result. $\Box$

It would be interesting to explore whether the techniques in this
note can be applied to evolving graphs \cite{b}, where a dynamic
notion of distance becomes relevant.

\section*{Acknowledgement}

The author is indebted to the referees for helpful comments. This
work is supported jointly by the National Natural Science Foundation
of China (11505127), the Shanghai Pujiang Program (15PJ1408300), the
Program for Young Excellent Talents in Tongji University
(2014KJ036), and the Fundamental Research Funds for the Central
Universities (0800219319).

\end{document}